\documentclass[12pt,letterpaper]{amsart}    

\usepackage[parfill]{parskip}    

\usepackage{graphicx}

\usepackage{amssymb} 

\usepackage{amsmath} 

\usepackage{amsthm} 

\usepackage{enumerate} 

\usepackage{fullpage}  

\usepackage{setspace} 

\usepackage{color}

\usepackage{wrapfig}

\theoremstyle{plain}

\begin{document}

\title{Distinguished graduates in mathematics of Jagiellonian University in the interwar period. Part II:  1926-1939}

\author[S. Domoradzki]{Stanis\l aw Domoradzki}
\address{Faculty of Mathematics and Natural Sciences,  University of Rzesz\'ow,  ul. Prof. S. Pigonia 1, 35-959 ,  Rzesz\'ow, Poland}
\email{domoradz@ur.edu.pl}
\author[M. Stawiska]{Ma{\l}gorzata Stawiska}
\address{Mathematical Reviews, 416 Fourth St., Ann Arbor, MI 48103, USA}
\email{stawiska@umich.edu}

\date{\today}                                           

\maketitle

\tableofcontents

\setcounter{tocdepth}{3}

\begin{abstract}
In this study, we present profiles of some distinguished graduates in mathematics of the Jagiellonian University from the years 1918-1939. We discuss their professional paths and scholarly achievements, instances of scientific collaboration, connections with other academic centers in Poland and worldwide, involvement in mathematical education and teacher training,  as well as their later roles in
Polish scientific and academic life.  We also try to understand in what way they were shaped by their studies and how much of Krak\'ow scientific traditions they continued.  We find  strong support for the claim that there was a distinct, diverse  and deep mathematical stream in Krak\'ow between the wars, rooted in classical disciplines like differential equations and geometry, but also open  to  new trends in mathematics. 
Part II concerns people who graduated after the university reform, getting the newly introduced master's degree in mathematics, in 1928-39.
\par W niniejszej pracy przedstawiamy sylwetki niekt\'orych wybitnych absolwent\'ow Uniwersytetu Jagiello\'nskiego w zakresie matematyki z lat 1918-1939. Omawiamy ich drogi zawodowe i osi\c agni\c ecia naukowe, przyk\l ady wsp\'o\l pracy naukowej, zwi\c azki z innymi o\'srodkami akademickimi w Polsce i na \'swiecie, zaanga\.zowanie w nauczanie matematyki i kszta\l cenie nauczycieli oraz ich p\'o\'zniejsze role w polskim \.zyciu  akademickim. Pr\'obujemy tak\.ze zrozumie\'c, w jaki spos\'ob zostali oni ukszta\l towani przez swoje studia  i na ile  kontynuowali krakowskie tradycje naukowe. Znajdujemy mocne dowody na poparcie tezy, \.ze w Krakowie mi\c edzywojennym  istnia\l \ wyra\'zny, zr\'o\.znicowany i g\l \c eboki nurt matematyczny, zakorzeniony w dyscyplinach klasycznych takich jak r\'ownania r\'o\.zniczkowe i geometria, ale r\'ownie\.z otwarty na nowe trendy w matematyce. 
Cz\c e\'s\'c II dotyczy  absolwent\'ow po reformie kszta\l cenia uniwersyteckiego, kt\'orzy uzyskali nowo wprowadzony stopie\'n magistra matematyki w latach 1928-1939.

\end{abstract}

\section{Introduction} After regaining independence in 1918, the Second Republic of Poland had to  build a unified modern state out of territories formerly under occupation of the three superpowers of Russia, Germany and Austro-Hungary. This required, among other things, creating a common educational system. Few academic schools existed continuously on partitoned territories;  Jagiellonian University was one of them. Some schools established  earlier were closed and then reviwed during World War I or after its end (University of Warsaw; Stefan Batory University in Vilnius). New institutions of higher educations were established, e.g., the Academy of Mining in Krak\'ow (1919). The first legal bill concerning the higher education was issued in  1920. More detailed regulations followed. On March 12, 1926, a decree of the Minister of Religious Denominations and Public Education was issued, which concerned the curriculum of studies and examinations in the field of mathematics for the master's degree.  The introduction of this degree was an innovation, at first optional for the students, but soon it  became an educational standard.  The ministerial decree stipulated that, during their course of studies, the students had to pass several exams (differential and integral calculus with introduction to analysis; analytic geometry; principles of higher algebra with elements of number theory; theoretical mechanics; experimantal physics; main principles of philosophical sciences; a block of two exams in pure or applied mathematics to be determined by the Faculty Council; and additionally one of a few subjects designated as ``auxiliary"). The final exam concerned general mathematical knowledge and was accompanied by the discussion of the master's thesis (\cite{Dy00}); for an in-depth discussion of the formation of the higher education system in Poland between the wars see \cite{Bar15}.

Andrzej Turowicz was the first at the Jagiellonian University to get the master's degree in mathematics:\\
\textit{``I was the first in Krak\'ow to get the master's diploma. When I enrolled [at the university], one could  pursue the old course, pre-master. I decided to do the master's degree. The second master [in mathematics] in Krak\'ow was [Stanis\l aw] Turski. He was two years younger than I.  [Zofia] Czarkowska (currently Mrs. Krygowska) was at the university along with me. I had a gifted classmate, Stefan Rosental, who however later became a physicist and finished his career as a vice-director of the Bohr Institute in Copenhagen."} (\cite{TuAU}\\

  In the years 1926-1939 lectures were given by S. Zaremba, A. Hoborski, A. Rosenblatt, W. Wilkosz, T. Wa\.zewski, S. Go\l \c ab, O. Nikodym, L. Chwistek,  J. Sp\l awa-Neyman, A. Ro\.za\'nski, F. Leja, J. Le\'sniak and S. K. Zaremba (\cite{Go64}). At the initiative of Wilkosz, two assistants (Jan Le\'sniak and Irena Wilkoszowa, cf. \cite{DS15a}) were employed and more lectures were enhanced with the recitation classes. An important part of the course of studies was teachers' training. Many students chose the teaching career. Because of the shortage of academic jobs  in Poland, a country with  nearly 35 million of population and about 40 fully accredited academic schools in 1938 \footnote{According to \cite{Wie02}, there were in total about 800 professorial positions and 2700 junior faculty positions in 1939; in mathematics, according to   \cite{Kur80}, there  were respectively 23 professorial chairs and 27 junior positions.}, even those who had talent and inclination for research started out as high school teachers, sometimes continuing for many years (e.g., A. Turowicz). A course in elementary mathematics from the higher standpoint was offered to address the needs of future teachers. It was taught by Jan Le\'sniak. \\

As recalled by Kazimierz Kuratowski (\cite{Kur80}) at the First Congress of Polish Science in 1951, an assessment was issued of the achievements of Polish mathematics in the inter-war period. It stated that  the greatest achievements were in functional analysis and topology; important contributions were made in real analysis, set theory and mathematical logic. Among other branches cited at the Congress were:

\textit{``Differential equations (in particular the results concerned with harmonic functions, the existence of integrals of partial differential equations of the second order, the qualitative theory of ordinary differential equations and the propersties of integrals of partial differential equations of the first order).\\
Geometry, together with transformation theory (in particular, the results concerning the invariants of surface bending, algebraic geometry, Finsler and Riemann spaces and the topology of geometric objects).\\
The theory of analytic functions (in particular the results concerning the approximation of functions by polynomials, the convergence of series of polynomials in many variables, univalent and multivalent functions)."}

These disciplines were precisely the strong points of mathematics at the Jagiellonian University, and moreover they were hardly represented anywhere else in Poland (although the report does not name any particular mathematical center in this context). Obviously these topics dominated not only the faculty's research, but also the students' master theses, although some topics in e.g., topology, measure theory or even functional analysis were represented, too (see the Appendix). In the period 1928-1939, over 135 people graduated form the Jagiellonian University with the master's degree in mathematics (see the Appendix for partial information; at the time of writing this article, we were not able to verify the data in full). Below, we present the profiles of those who made their mark on Polish scientific and academic life.

\tableofcontents

\setcounter{tocdepth}{3}

\section{Profiles}

\subsection{Stanis\l aw Krystyn Zaremba (1903-1990)} 

\begin{wrapfigure}{L}{0.45\textwidth}

\begin{center}
\includegraphics[width=0.45\textwidth]{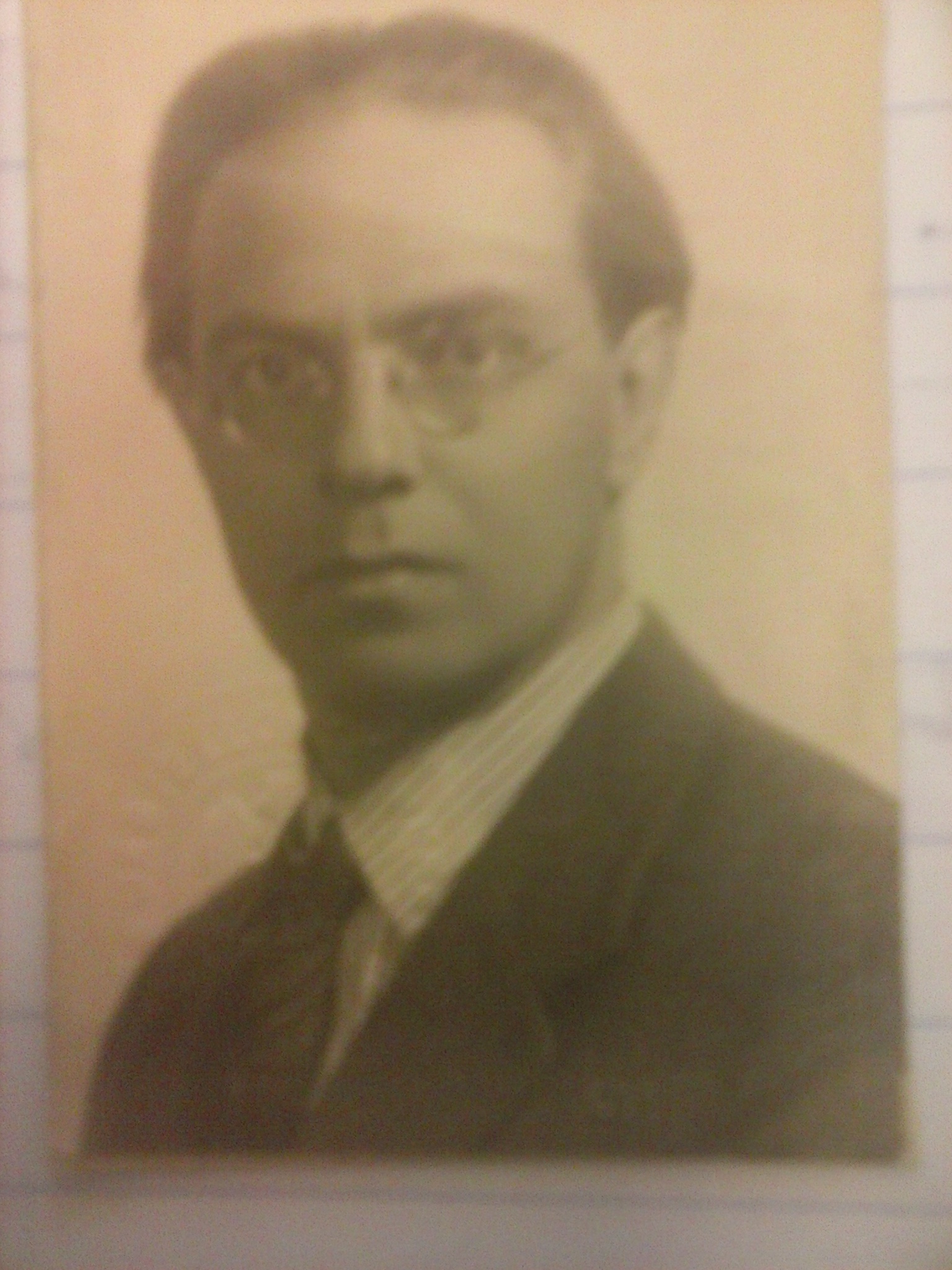}
\end{center}

\end{wrapfigure}

Born in Krak\'ow, a son of mathematician Stanis\l aw Zaremba (a professor of Jagiellonian University) and a Proven\c cal woman Henrietta Leontyna ne\'e Cauvin. After finishing high school with science-oriented curriculum in 1921, he started studying mathematics at Jagiellonian University. Following in his father's footsteps, he continued his studies at the Sorbonne in Paris in the years 1924-27. Because of health problems he returned to Krak\'ow, where he got  master's degree in mathematics from Jagiellonian University in 1929. He edited lectures of Professor Jan Sleszy\'nski, which were later published as a two-volume \emph{Proof Theory} (in 1923 and 1929). Since 1929 he was an assistant at the Stefan Batory University in Wilno (Vilnius). There he got PhD degree on the basis of the thesis  (\cite{Za39}) ``Sur l'allure des int\'egrales d'une \'equation diff\'erentielle  ordinaire du premier ordre dans le voisinage de l'int\'egrale singuli\`ere" (supervised by Juliusz Rudnicki).  He also mentored a distinguished student Duwid Wajnsztejn, who went on to obtain PhD in Krak\'ow. In 1936 Zaremba got his habilitation at Jagiellonian University on the basis of the thesis ``On paratingent equations" (\cite{Za35}, \cite{Za36}), in which he introduced paratingent equations, a generalization of differential equations nowadays known as differential inclusions. About the same time similar relations were independently studied by Andr\'e Marchaud  (\cite{Gor11}). This generalization later allowed Tadeusz Wa\.zewski and others to build natural foundations of optimal control theory. Since 1937 he was back in Krak\'ow, first as an \emph{adiunkt}, later as a \emph{docent} at  Jagiellonian University.\\

When World War II broke out, Zaremba returned to Vilnius, then under Lithuanian control. After Lithuania was annexed by the Soviet Union in 1940, he went to Stalinabad (now Dushanbe, in Tajikistan), where he worked as a professor of mathematics in Pedagogical Institute. Along with the Polish Army (formed from Polish nationals in the USSR under the command of General W\l adys\l aw Anders) he went first to Persia, then to Palestine, where he taught in high schools for the army. In 1946 he worked at the University of Beirut. Fearing persecution or even death from the new   communist Polish authorities, he decided to stay in the West. Since 1946 he lived in Great Britain. Until 1952 he was a professor at the Polish University College in London. Then he became a mathematical consultant for Boulton Paul Aircraft Ltd. in Wolverhampton. This position-- as he acknowledged himself-- gave him an opportunity to familiarize himself with the theory of stochastic processes and start research on this subject. In 1950s  he collaborated with Zbigniew \L omnicki, a graduate in mathematics and physics of the Lw\'ow University and a fellow emigr\'e, on the theory of time series. They published 8 joint papers. (\cite{Da98}). In 1954 he took part in the International Congress of Mathematicians in Amsterdam, where he gave a 15-minute talk ``Spacing problems in Abelian groups" on pioneering application of group theory to communication theory (\cite{ICM}). In the years 1958-69 (with a yearlong break in 1966/67) he lectured at the University of Wales.  He spent the years 1969-76 in North America (Madison, WI, and Montr\'eal, Qu\'ebec). In 1976 he returned to Wales. Before the martial law in 1981 he frequently visited Poland, in particular Warsaw and Krak\'ow.  He died in Aberystwyth. (See also \cite{Du12}, \cite{PB03Za}, \cite{ZaXX}.)\\ 

In the years 1925-37 Zaremba was an active mountaineer. He made many first routes and winter ascents as a pioneer of snow mountaineering in Tatra mountains. He  was on boards of mountaineering societies. He published literary accounts of his expeditions in the journals ``Wierchy", ``Krzesanica" and ``Taternik" (of which he was an editor in 1929-30). While in emigration, he climbed e.g. in the Hindu Kush range. At his wish, his ashes were scattered over Tatra and the mountains of Wales (\cite{RPP73}).\\

In recent years, there has been substantial interest in so-called Zaremba's conjecture (\cite{Kon13}). Stated in \cite{Za71}, it was motivated by his search for lattice points suitable for quasi-Monte Carlo methods in numerical integration and  postulates the following: there is an universal integer constant $K > 0$ such that for every integer $d > 0$ there exists an integer $b$ coprime with $d$, $1 \leq b < d$, such that all partial quotients of the continued fraction expansion $b/d=[0;a_1,a_2,...,a_k]$ satisfy $a_i \leq K$. Zaremba also conjectured that $K=5$. In the paper \cite{BK14} the problem was reinterpreted in terms of properties of the orbit of the vector $(0,1)$ under some semigroup of matrices and the conjecture was proved for almost all $d$ (in the sense of density) with $K=50$. \\

\subsection{Andrzej Turowicz (Fr. Bernard OSB; 1904-1989)} Born in Przeworsk, in the family of August, a judge, and Klotylda ne\'e Turnau. Initially homeschooled, he finished King Jan Sobieski Gymnasium in Krak\'ow. After his \emph{matura} he studied mathematics at Jagiellonian University, in the years 1922-28. He was the first graduate to obtain the degree of master of philosophy in the area of mathematics. In 1931 he took a high school teacher qualifying exam and started teaching in  schools in Krak\'ow and Mielec. While in Krak\'ow, he combined teaching school with academic activities. In the years 1929-30 he was an assistant in the Chair of Mathematics of the Academy of Mining in Krak\'ow, substituting for Stanis\l aw Go\l \c ab, who went to Netherlands on a scholarship.  In 1937 he got a position of the senior assistant in the Chair of mathematics of Lw\'ow Polytechnics, where he worked until 1939. Here is how he recalled the circumstances (\cite{TuAU}, cassette 1b) of his appointment:\\ 


\begin{wrapfigure}{L}{0.45\textwidth}
\centering
\includegraphics[width=0.45\textwidth]{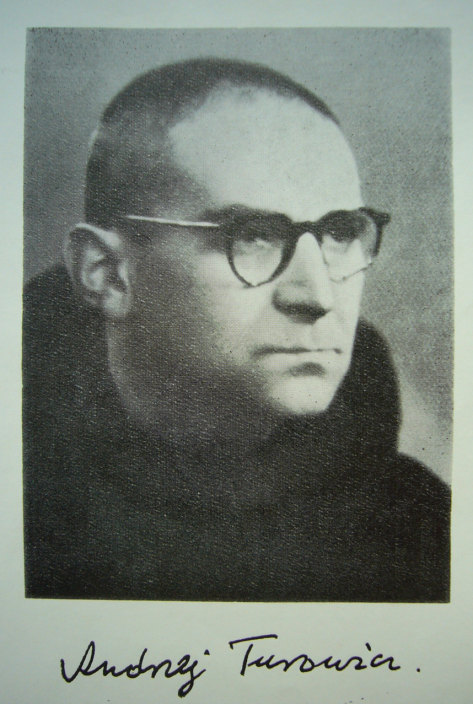}

\end{wrapfigure}

\textit{``After the J\c edrzejewicz Brothers  reform, [Antoni] \L omnicki submitted a geometry textbook for high schools. The Ministry gave me this text for refereeing, [as] I taught high school in Krak\'ow. I did a very detailed report; I went over all problems, I wrote which ones were too difficult. The authors were anonymous, [but] I guessed it was \L omnicki. I remembered his geometry and trigonometry [textbooks] from Austrian times. However, the referees' names were made known to authors.  When he read my report,  [\L omnicki] sent me an offer to become an adiunkt. He wanted to bring me where he was for this report."}\\

After the Soviet-style reorganization of the Polytechnics, with Ukrainian as the language of instruction, Turowicz  taught mathematics at the Faculty of Architecture, then at the Faculty of Mechanics. In 1941 he returned to Krak\'ow, where he worked as a clerk in the Chamber of Industry and Commerce until 1945. At the same time, he taught analytic geometry at the underground university and took part in clandestine sessions of the Polish Mathematical Society.\\

The day after the entrance of the Red Army to Krak\'ow, January 18, 1945, Turowicz crossed the frozen Vistula river to the Tyniec Abbey. Ten days later he entered the Benedictine order, taking the name of Bernard. In the years 1946-50 he studied theology. He was ordained a priest in 1949. In 1946 Turowicz obtained at the Jagiellonian University the PhD degree on the basis of the thesis ``On continuous and multiplicative functionals" (published as \cite{Tu48}).  Tadeusz Wa\.zewski, the supervisor,  noted in his report that Turowicz's thesis was an evidence ``of deep mathematical culture and revealed rare philosophical sense of its author in treating a problem." The  thesis answered a question posed by Stefan Banach and Meier Eidelheit. Here is how Turowicz related the story about his work on the problem:

\textit{I had an incident with Banach like this:  for  a meeting of the mathematical society, I proposed a talk on multiplicative and continuous functionals (my proposal was in spring 1939; I finished the work after the war). I am delivering the talk and Banach enters the room, slightly late, with an incredibly sullen face. I noticed that  Banach was angry. He listened with extreme attention and his face changed. When I finished, Banach took the floor and said: ``I also dealt with this problem; you did it in a totally different way, and you did it well." I received his opinion with gladness. The next day after this meeting, Sto\.zek (who was not at the meeting) asked: ``Was Banach there?" [I said]  ``Yes, he was." [He said] ``I did not want to scare you in advance; he was very angry when he found out what you were to talk about. He said: `I am dealing with this; I must have told someone, and [now] Turowicz is presenting it as his own.'" Banach came with the intention of giving me a hard time. Luckily the idea of the proof was completely different [from his], therefore [he] praised me [and] did not make a scene. (...) Since then, Banach was very friendly towards me.}
(\cite{TuAU}, cassette 2b)

In the years 1946-52 and 1956-61 Turowicz gave lectures in mathematics at the Jagiellonian University, initially as the replacement for W\l adys\l aw Nikliborc, who moved to Warsaw before the end of the term. In 1949 he lectured on algebra to the first-year students. Czes\l aw Olech,  who was taking the class, had the following memory of Turowicz (\cite{Ol10}): \textit{``He commuted for our lectures from the monastery in Tyniec near Krak\'ow, where he resided. He almost dashed into the lecture room, made the sign of the cross and filled the blackboard with legible text. If someone could take exact notes of it-- and there were women in the class who could-- then we had a `textbook' for the exam."} The lecture notes from the years 1946-48 were indeed published (internally) as ``The theory of determinants and matrices with applications to the theory of linear equations and the forms of 1st and 2nd degree" (1949), presenting some contents for the first time in the Polish mathematical literature.\\

At the beginning the new communist authorities did not object to Turowicz's  appointment, as they took efforts  to rebuild academic life in Poland facing the  shortage of qualified scholars. (Later, Turowicz gave credit to his classmate Stanis\l aw Turski, who worked for the Ministry of Education, for  signing an appropriate permission.)  However, in the years 1952-56-- the Stalinist period in Poland-- there was no place for a priest at a state institution, so  he taught mathematics for philosophers at the Catholic University of Lublin. 
In 1954 the institution applied for granting him the title of \emph{docent} (i.e., an independent scientific worker) on the basis of his scholarly output. It was a legitimate procedure at that time. At the request 
of the dean of the Faculty of Philosophy evaluations were written by Hugo Steinhaus and Tadeusz Wa\.zewski. Steinhaus wrote, among other things, that, when Turowicz was in Lw\'ow, \textit{``(...) I, along with other Lw\'ow mathematicians, had an impression that we dealt with a young mathematician who would develop his creative abilities in the right conditions."} Mentioning Turowicz's publications, he wrote that \textit{``(...) all these works are an evidence of mathematical cultivation and creative capabilities of the author. (...) I add that all colleagues whose opinion I asked said without reservation that Dr. Andrzej Turowicz fully deserves the title of docent. I must also repeat here a general opinion about great personal qualities of Dr. Turowicz, who enjoys universal respect in the circles of his acquaintances and colleagues."}  Despite very good evaluations, the application took a long time and was ultimately denied in 1957.  At Wa\.zewski's insistence, in 1961 Turowicz obtained a position in the Mathematical Institute of the Polish Academy of Sciences. It was a research position; the employees of the Institute did not have direct contact with students. In 1963 he got his habilitation on the basis of a series of papers about ``orientor fields" (i.e., differential inclusions) and their applications to control systems. He became an extraodinary professor in 1969. In the years 1970-73 he taught in the doctoral study program organized by the Faculty of Electrotechnology, Automated Control  and Electronics of the Academy of Mining and Metalurgy. The lecture notes for some classes he gave there were published in a book form under the title ``Matrix Theory" (\cite{Tu85}). At the Tyniec Abbey, for a few years he  taught history of monasticism to  candidates for the holy orders. He retired from  academic positions in 1974.\\

Turowicz's scholarly output is very diverse and spans functional analysis, differential equations, control theory, probability, linear algebra, logic, game theory, convex geometry, algebra, functions of one complex variable  and numerical analysis.   In Lw\'ow he wrote one joint paper with Stefan Kaczmarz (\cite{KaTu}). He also collaborated with Stanis\l aw Mazur, but their joint results were never published, even though Mazur found the manuscript after World War II. The work concerned a generalization of Weierstrass' theorem on approximation of continuous functions by polynomials, akin to what is now known as the Stone-Weierstrass theorem (Stone proved his versions in 1937 and 1948). The reasons for not submitting the paper for publication were twofold, and clear to those who knew Mazur, including Turowicz himself (\cite{CP88}, \cite{Tu95}). First, Mazur was always striving for the best possible version of his results, and delayed submissions in hope of improving them. Second, he was a committed communist, so he distanced himself from a former colleague who became a priest. In the written evaluation of Turowicz's output in 1950s  Wa\.zewski mentioned his joint results with Mazur and expressed regret that they were unpublished. Turowicz successfully collaborated with other scholars, e.g., with H. G\'orecki  on applications of mathematics to automated control theory. They published several joint papers and a monograph ``Optimal Control" (\cite{GoTu70}). Another monograph by Turowicz, ``Geometry of Zeros of Polynomials" (\cite{Tu67}) published in 1967, concerns polynomials in one complex variable and discusses the number of real zeros, localization of zeros of polynomials and their derivatives and other related topics. It was written with engineering applications in mind, but it contains many classic mathematical results.\\

Turowicz was active in the Polish Mathematical Society, starting from 1927. In the years 1973-75 he was the president of the Krak\'ow branch of the Society. Since 1978 he collaborated with the Committee for History of Mathematics by the General Management of the Polish Mathematical Society. He had an incredible memory and a gift of storytelling, without taking himself or the surrounding world too seriously. He was willing to meet with students and give interviews. He was also regarded as a moral authority. (See also \cite{Du12}, \cite{PB03Tu}.)\\

\subsection{Stanis\l aw Turski (1906-1986)} \begin{wrapfigure}{L}{0.45\textwidth}

\begin{center}
\includegraphics[width=0.45\textwidth]{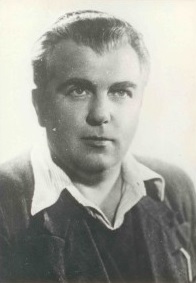}
\end{center}

\end{wrapfigure}

Born in Sosnowiec,  he finished high school there. In the years 1924-28 he studied physics and astronomy at the Jagiellonian University.  During his studies he worked as a schoolteacher. His diploma thesis ``A new method of determining precession coefficients" was awarded a prize by the minister of education. Since 1927 he worked as a mathematician,  starting at the level of an assistant, at the Jagiellonian University and at the Academy of Mining and Metalurgy.  He also gave lectures in mathematics in Krak\'ow Pedagogium. He obtained his PhD under the supervision of Witold Wilkosz, presenting a thesis ``On a generalization of theorems on uniformity of integrals of a hyperbolic equation". Arrested in 1939 in Sonderaktion Krakau, he was imprisoned in the concentration camps of Sachsenhausen and Dachau. After his release in 1941 he took part in clandestine teaching at an academic level. \\

After the World War II, Turski engaged in rebuilding academic life in Poland, as a supporter of the communist party and its program. Nominated by  Minister Stanis\l aw Skrzeszewski  (a recipient of PhD in logic from Jagiellonian University), he lead a group dispatched by the Ministry of Education to  reactivate the    Gda\'nsk University of Technology in 1945  as a Polish-staffed  institution replacing a  German academic-level polytechnic school. He was an extraordinary professor of mathematics and a rector (president) in the years 1946-49  (\cite{Wit}). He also became a parliament member in 1947. In 1949 he was called to work in Warsaw, at the  University of Warsaw and in the Ministry of Education. He became an ordinary professor of mathematics in 1951 and got habilitation in 1953. In 1954 he took part in the International Congress of Mathematicians in Amsterdam as a delegate of the Polish Academy of Sciences. In the years 1952-69 he served as the rector of the University of Warsaw. During his term, in March  1968, student protests against communist authorities erupted. As a result,   34 students were expelled and 11 were suspended   from the university,  and  the professors were officially prohibited from  participating in  students' rallies. Other repressions followed.\\

 Turski's work before World War II concerned partial differential equations, number theory and functions of a complex variable. Out of his publications of that period,  three were joint with Alfred Rosenblatt, a \emph{docent} in Krak\'ow (\cite{RT35}, \cite{RT36b}, \cite{RT36c}).  After the war he published a few papers applying mathematical methods in mechanics of solids. His paper with Jerzy Nowi\'nski \cite{NT55} contained a numerical solution to a system of differential equations occurring in  elasticity theory obtained with the use of ARR (Differential Equations Analyzer)-- the first analog computer constructed in Poland in 1953. In 1963 Turski  arranged for an exhibition, followed by purchase in 1964,  of an ALGOL-running computer from Denmark  and for staff-training courses, which lead to creation of the  Unit of Numerical Computations at the University of Warsaw.  To reflect the emergence of a new direction in research and education,  his Chair of General Mathematics was renamed the Chair of Numerical Methods. These two institutions were later combined to give rise to the Institute of Computer Science. (\cite{MS00}). Turski retired in   1976.\\

\subsection{Antoni Nykli\'nski (1906-1964)} 

\begin{wrapfigure}{l}{0.3\textwidth}
\begin{center}
\includegraphics[width=0.3\textwidth]{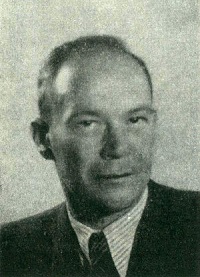}
\end{center}
\end{wrapfigure}
 Born in Krak\'ow,  he studied mathematics at the Jagiellonian University in the years 1925-1930, obtaining the master's degree in 1932. Afterwards he worked in high schools and collaborated with the Chair of Mathematics of the Academy of Mining on research in differential geometry. He translated into Polish a popular book on mathematics by Egmont Colerus (\cite{Col38}); the preface to the Polish edition was written by Stefan Banach. In 1939 Nykli\'nski was arrested, imprisoned at the Montelupi Street, then in Wi\'snicz Nowy, and then taken to  the concentration camp in Auschwitz. After the World War II he resumed his work in high and academic schools. He conducted lectures in mathematics at the Preparatory Study of the Jagiellonian University since 1951. In 1956  he was nominated to the post of adiunkt in the Chair of Mathematics at the Faculty of Electrification of Mining and Metallurgy of the Academy of Mining and Metallurgy. In 1962 he got his PhD degree at the Faculty of Finance and Statistics of the Main School of  Planning and Statistics in Warsaw. His research interests and educational activities focused on linear programming and probability. He died in Krak\'ow. (\cite{ZyAGH64})

\subsection{Czes\l aw Kluczny (1908-1979)} Born in Strzemieszyce Wielkie, he finished a gymnasium in Olkusz in 1927. In the years 1927-32 he studied mathematics at Jagiellonian University,   obtaining a master's degree. In 1932 he started working as a high school teacher in Radom. In 1942 he was arrested by Gestapo and sent to concentration camps-- first  Auschwitz-Birkenau, then Mauthausen. He returned to Poland in 1945 in very poor health, recovering for a year. In 1946-50 he was employed by the Silesian Technical Scientific Institutions in Katowice. In 1950 he started working for the Gliwice Polytechnic, where he remained until his retirement in 1976.  He also lectured at Silesian University in Katowice and at the Pedagogical College in Cz\c estochowa. In 1959 he obtained PhD degree in mathematics at Jagiellonian University under the supervision of T. Wa\.zewski. In 1961 he got habilitation at the Maria Curie-Sk\l odowska University in Lublin. He became an extraordinary professor in 1971.\\

Kluczny worked in qualitative theory of ordinary differential equations (\cite{Klu60}, \cite{Klu61}, \cite{Klu62}). He  supervised 6 PhD degrees (3 of which were interrupted because of his death) and 1 habilitation. He was a co-founder of the Gliwice (later Upper Silesian) branch of the Polish Mathematical Society. \\

\subsection{W\l adys\l aw Benedykt Hetper (1909-1941?)} 
\begin{wrapfigure}{l}{0.25\textwidth}
\begin{center}
\includegraphics[width=0.25\textwidth]{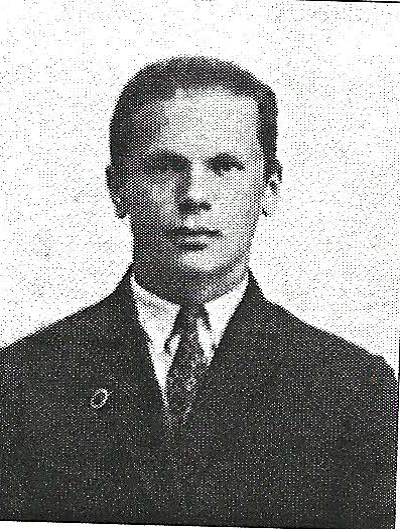}
\end{center}
\end{wrapfigure} 
Born in Krak\'ow, he finished the King Jan Sobieski 3rd Gymnasium there. Then he studied  mathematics at Jagiellonian University in 1927-32, obtaining master's degree.  His master thesis concerned integral equations, but soon he became interested in logic under the influence of Leon Chwistek. Along with Jan Herzberg and Jan Skar\.ze\'nski, he collaborated with Chwistek on  his program of establishing consistent foundations of mathematics (\cite{CHH33a}, \cite{CHH33b}, \cite{CH38}).\footnote{The following  description of Chwistek's program can be found in \cite{Jo63}; see also \cite{Wol15}: ``The first step was the creation of what Chwistek called ‘elementary semantics’, which, besides its name, has nothing in common with semantics in Tarski's sense. Chwistek claimed that in Hilbert's metamathematics there is contained intuitive semantics, i.e. the rules for the construction of the simplest possible expressions from given elements (letters or signs). This intuitive semantics is formalised and expanded into a system of syntax in terms of which the propositional calculus and the theory of classes are constructed. On this basis the axiomatisation of classical mathematics which assumes no non-constructive objects is finally undertaken. If successful, and this matter must be left to the mathematician to judge, it would provide a proof of the consistency of mathematics. In this, more than in anything else, lies the importance of Chwistek's system."} 
 When   Chwistek took the Chair of Logic of the Jan Kazimierz University in 1933,  his students followed him there.  Hetper was the most active of them, publishing several papers, in which he paid attention to the latest developments in logic. 
E.g., in \cite{He38}, he proposed structural rules (in the style of sequent calculus of Gerhard Gentzen, introduced in 1934) for a propositional calculus written in so-called Polish notation, introduced  by Jan \L ukasiewicz in 1924. He proved consistency and completeness of his system, comparing his methods to those of Hilbert. In the introduction he credited Witold Wilkosz, who apparently worked on similar problems.\footnote{``It has been known  to me through private communication that similar problems were studied by Dr. W. Wilkosz, professor of Jagiellonian University; however, I do not know of any of his publications or results in this direction." (\cite{He38}).}  While working on his PhD (supported by a government scholarship), Hetper shared a room with Mark Kac, with whom he became friends.  They had intellectual discussions, played chess and went cross-country skiing together. A devout Catholic, Hetper  represented a positive example of Christianity to his secular Jewish friend. Hetper and Kac  got their PhD degrees and had them conferred in a double ceremony on June 5, 1937 (\cite{Kac85}.  Soon Kac left Poland and Hetper went on to get his  habilitation in 1939 (the thesis \cite{He38a} was printed in 1938). 
When the World War II broke out, he fought in the September campaign in 1939 as an ensign of infantry reserves. He escaped from German captivity, but was  arrested by Soviets at an  attempt of illegal border crossing.   Because  he carried with him a mathematical manuscript, he was accused of espionage  and  imprisoned.  (\cite{Pio84}). His last known address was the Starobielsk  camp and his last letter to his family was dated 1941. (We thank Professor Roman Duda for this information). Hetper died probably in 1941.\\

\subsection{Danuta Gierulanka (1909-1995)} \begin{wrapfigure}{L}{0.45\textwidth}

\begin{center}
\includegraphics[width=0.45\textwidth]{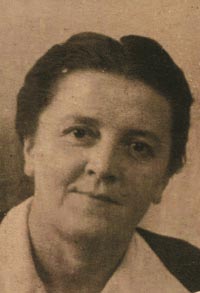}
\end{center}

\end{wrapfigure} Born on June 30, 1909, in Krak\'ow, in the family of Zofia ne\'e Romanowska and Kazimierz Gierula, a civil servant in the Ministry of Communication until 1939. In the years 1927-1932 she studied mathematics in the Faculty of Philosophy at the Jagiellonian University. She obtained a master's diploma  in philosophy in the field of mathematics for the thesis ``Periodic solutions of differential equations." After graduation she enrolled again in the University, in its Pedagogical Study, in order to prepare herself properly for a teaching licence examination, which she passed in 1933. She taught mathematics, physics, chemistry and propaedeutics of philosophy in gymnasia: Landowners' Gymnasium of Benedictine Nuns in Stani\c atki near Krak\'ow, Humanistic Gymnasium of Mary's Institute  in Krak\'ow. In 1938 she started doing research in psychology, under the supervision of W\l adys\l aw Heinrich, on the problems of psychology of thinking, and more precisely, on forming geometrical notions. During World War II she supported herself by giving private lessons and  working in  offices of commerce. She took part in clandestine teaching along with her younger brother Jerzy, later a renowned physicist.\\

In 1945 Gierulanka was nominated to the post of the senior assistant in the Laboratory  of Experimental Psychology of the Jagiellonian University. She combined these duties with teaching mathematics in one class in H. Ko\l \l \c ataj Lyceum in Krak\'ow. She initiated research on thought processes in mathematics. In 1947 she obtained a doctorate at the Humanistic Faculty of Jagiellonian University for the thesis ``On acquiring geometrical notions" (\cite{Gie58}).(She published a book based on her thesis in 1958). The primary subject of the examination was psychology, the secondary one was mathematics. The examiners were professors W. Heinrich, Stefan Szuman and Tadeusz Wa\.zewski (two psychologists and a mathematician). The mathematician's questions were interesting and showed his extremely favorable attitude to the problems of instruction: ``Non-euclidean geometry", ``On criteria of good elaboration of mathematical theorems", ``Characteristics of mathematical talent". Stefan Szuman wrote in his report: \textit{``She investigated the psychological process of forming clear and rigorous notions in geometry by the pupils."} Gierulanka passed the doctoral examination with  distinction. In 1953 she was transferred to the post of adiunkt in the Chair of Mathematical Analysis in the Faculty of Mathematics, Physics and Chemistry at the Jagiellonian University. It was suggested that she would prepare there a candidate's thesis (then an equivalent of a PhD thesis) in mathematical analysis incorporating psychological research, which would be a   counterpart of her PhD thesis regarding university teaching. The plan failed. In 1957 Gierulanka returned to the Laboratory  of Experimental  Psychology for a year. In 1958 she became an adiunkt in the Chair of Philosophy. She obtained habilitation in 1962 on the basis of the thesis ``The problem of specificity of mathematical cognition" (\cite{Gie62}), but she did not get the position of a docent in the  Chair of Philosophy. She was transferred to the Chair of Psychology, where she worked until her retirement in 1971. She retired with the sense of injustice:  the program of doctoral studies in university teaching for assistants in various academic disciplines was not launched and the university authorities did not show recognition of her work and achievements.  She died in Krak\'ow on April 29,  1995. \\

Gierulanka's scientific path can be best described in her own words. In the ``Information on my previous scientific work" attached to the application for habilitation she wrote:\\

\textit{``Influenced by lectures and seminars conducted by Prof. R. Ingarden, in which I participated regularly since 1946, I was getting an even broader view of the philosophical problematics related to my psychological work. A problem in which I have taken stronger and stronger interest since 1948 is the problem of specificity of mathematics; in its solution I would see a natural complement of the work on acquiring geometrical notions. [This work] gave an idea-- thanks to investigation of the course of suitable psychical thought processes-- only about psychological sources of the paradoxical opposition between the fundamental clarity and comprehensibility of mathematics and the actual state of its being comprehended. To explain it fully it is necessary to realize what the specific character of mathematics and of the cognitive means it employs consists in. Seeing Cartesian mental intuition and what is called clara et distincta perceptio as a cognitive activity  typical for mathematics, I analyzed this notion. (...) Another kind of cognitive processes very relevant for mathematical cognition are processes of understanding. I have been concerned with the problems of understanding  since 1951, conducting for 2 years research in Laboratory of Experimental Psychology concerning primarily understanding of texts. With the subsidy from Scientific Pedagogical Society I conducted research parallel to this, concerning learning mathematics from textbooks, including analysis and criticism of school textbooks in geometry in use at that time. (...) Because of the financial difficulties of the Scientific Pedagogical Society the research ceased, and the partial results obtained were not published; those concerning textbooks became obsolete when change occurred."}. The problematics undertaken by Gierulanka was very comprehensive. Later on, these issues were dominated by cognitivist psychology. The book \cite{Zaj12} is an attempt at re-reading Gierulanka's work.  \\

In her habilitation, Gierulanka addressed the problem of mathematical cognition as a philosopher, although she also used some previously collected psychological materials. The referees of her scientific output were Professors Zofia Krygowska, Izydora D\c ambska, Tadeusz Cze\.zowski and Roman Ingarden.  Gierulanka analysed mathematical perception and deduction from the phenomenological standpoint and  described attempts at systematization and unification  of mathematics. She did not consider reduction to set theory as true solution of the problem of unification of mathematics. She criticized Bourbakist mathematics, blaming it  for, among other things,  being arbitrary in constructing systems of axioms and  making unnatural generalizations.   However, she saw some possibilities for applying the notion of mathematical structure (\cite{Biel02}).  Krygowska considered her presentation of the  state of contemporary mathematics  exaggerated and tried to defend the Bourbakist approach by pointing out that it can reflect an actual course of mathematical creative processes as reported by some distinguished mathematicians. Ingarden agreed with Gierulanka's reference to Descartes, considering his epistemology still relevant for contemporary mathematics. He wrote:\\

\textit{``One needs to realize that mathematics is not only a certain set of theorems and methods used over  the last decades, but that it is a certain historical creation, evolving over at least last three centuries, precisely since the Cartesian reform. During that time not only did mathematics encompass even more new domains of study, but it also significantly kept changing its  methods, the understanding of its role and of its ultimate sense, undergoing a series of internal crises (e.g., emergence of non-euclidean geometries, antinomies at the end of 19th century, and finally G\"odel's theorems in 1930s) as well as a series of external crises through attacks of various forms of modern European scepticism, e.g., Hume's attack, various forms of positivism and empiricism of 19th century, up to Vienna neo-positivists in 20th century, who made mathematics a system of tautologies. Because of this the problem of specificity of mathematical cognition is extremely complicated and does not allow one to restrict considerations of mathematics  only to the form it has had in the latest decades. It cannot be excluded that the cognitive tendencies represented by Descartes-- despite differing current  views-- do not lose their relevance."}\\ 

Gierulanka was also an active editor and translator. She took part in translating some of Ingarden's works from German to Polish and in editing his collected works. She translated the first and second volume of Husserl's ``Ideas of pure phenomenology and phenomenological philosophy" and (together with Jerzy Gierula) ``On the problem of empathy" by Edith Stein.

\subsection{Adam Bielecki (1910-2003)}  \begin{wrapfigure}{L}{0.45\textwidth}

\begin{center}
\includegraphics[width=0.45\textwidth]{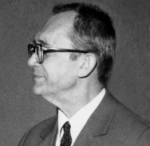}
\end{center}

\end{wrapfigure}

Born on February 13, 1910, in Borys\l aw (Drohobycz county, Lvov voivodship). In 1928 he finished Hoene-Wro\'nski Gymnasium in Krak\'ow. Simultaneously with the high school course he studied piano and theory of music at the Krak\'ow Conservatory. He enrolled in the Jagiellonian University to study mathematics. As a student, he taught acoustic in a private \.Zele\'nski School of Music in Krak\'ow.  He finished his studies in 1931, obtaining the master's degree in philosophy in the field of mathematics. In 1935 he obtained the doctorate on the basis of the thesis ``On integral representation of $m$-dimensional surfaces contained in the $n$-dimensional euclidean space by implicit functions" (\cite{Bie35}). The problem-- of equivalence of the implicit and parametric representations-- was posed  by Witold Wilkosz, who supervised the thesis. Bielecki learned about it from Tadeusz Wa\.zewski and solved it using the technique of $\mathcal{C}^\infty$-partition of unity, later reintroduced and refined by Laurent Schwartz  (\cite{KKZ04}).  In the years 1935-36 he worked in the Theoretical Physics Seminar at the Jagiellonian University, becoming a senior assistant in the Chair of Theoretical Physics in 1936. He collaborated with his colleagues, publishing 1 joint paper with Stanis\l aw Krystyn Zaremba (\cite{BieZa36}), 1 with Jan Weyssenhoff  (a theoretical physicist; \cite{BW36}) and 1 with Weyssenhoff and Myron  Mathisson (\cite{BMW39}). He also published two volumes of poetry. He was arrested in the Sonderaktion Krakau and taken first to the Sachsenhausen-Oranienburg concentration camp, then to Dachau. Released in April 1940, he returned to Krak\'ow. He supported himself first by giving private lessons, then-- from September 1942 to January 1945-- by part-time teaching at the Vocational School of Construction in Krak\'ow. At the same time, starting in 1942, he became active in organizing clandestine teaching in the underground Jagiellonian University, holding some classes in his private apartment, preparing lecture notes for students and lending them books from the remainder of the Library of the Laboratory of Theoretical Physics. He was also a member of the underground research group in theoretical physics led by J. Weyssenhoff.\\

	After the war, Bielecki decided to ``give his strengths to the Mathematical Institute" (\cite{Go64}). In 1945, he worked first as a senior assistant, then as an adiunkt of the I Mathematical Laboratory at the Jagiellonian University. From 1945 to 1947 he was a deputy professor and a head of the Chair of Mathematics at the Faculty of Engineering of the Academy of Mining in Krak\'ow. He was also strongly involved in the activities of the Polish Mathematical Society, to which he belonged since 1931 (taking part in clandestine scientific meetings during the occupation). In 1947 he was called to Lublin and assumed (as a deputy professor) the Chair of Mathematical Logic and Foundations of Mathematics  at the Faculty of Mathematics and Sciences of the Maria Curie-Sk\l odowska University. In 1949 he got habilitation on the basis of the thesis  concerning differential equations and differential inclusions, \cite{Bie48}. In 1959, after the death of Mieczys\l aw Biernacki, Bielecki took over as the head of the Collective Chair of Mathematics. In order to save the mathematics program from liquidation, he supervised a few PhD theses and supported 3 habilitations in a 3-year period of time. Some of the candidates started their research under the direction of Biernacki and worked in the theory of univalent functions in one complex variable. Not only did they finish their theses, but Bielecki was able to adapt his interests  in a way that allowed him   to write joint papers with them, concerning mainly subordination theory (e.g., \cite{BieLe62}). This is still an active research topic, especially in relation with the Loewner differential equation (which was used in 1985  in Louis de Branges's proof of the Bieberbach conjecture and whose stochastic version, introduced  by Oded Schramm in 2000 and later studied by him together with Gregory Lawler and the Fields medalist Wendelin Werner, found applications in statistical mechanics and conformal field theory). Besides university teaching and supervising PhD candidates (11 over the whole career), Bielecki was active in curriculum development and teachers' education. He organized post-graduate courses for mathematics teachers and qualifying exams for those teachers who did not have a master's degree in mathematics. In 1970s  he delivered  lectures  and created lesson plans for teachers as a part of Radio and Television Teachers University educational program. His presentations were later followed by articles containing in-depth discussion of mathematical and educational issues, published in the biweekly ``O\'swiata   i Wychowanie". He retired in 1980, but continued part-time teaching until 1991. He died in Lublin on June 10, 2003.\\

Bielecki's main mathematical interest were differential equations and differential inclusions. His best-known and most-cited result (\cite{Co94}) is a method of proving theorems on existence of solutions of differential and integral equations. The method consists  in a suitable change of a norm in the relevant space of functions so that a certain operator becomes a contraction and Banach Fixed Point Theorem can be applied (\cite{Bie56a}, \cite{Bie56b}).  In another paper (\cite{Bie57a}) Bielecki extended Wa\.zewski Retract Theorem to differential inclusions. He also worked on various aspects of geometry, publishing, among others,  the paper \cite{BMW39} motivated by aplications to the general theory of relativity and  a joint paper with S. Go\l \c ab (\cite{BieGo45}) concerning characterization of Riemannian space among Finsler spaces by the properties of the angular metric. His notable result in foundations of geometry (\cite{Bie56c}, \cite{Bie57b}) concerns reducing  the number of axioms given by Hilbert for Euclidean geometry (while weakening some of them) and the proof of independence of the system obtained.\\

\subsection{Antoni Bielak (1910-1991)} \begin{wrapfigure}{l}{0.3\textwidth}

\begin{center}
\includegraphics[width=0.25\textwidth]{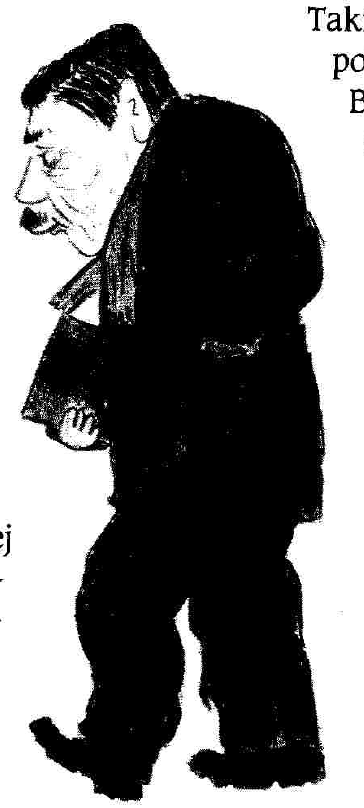}
\end{center}

\end{wrapfigure} Born on June 10, 1910, in Krak\'ow, in the family of Antoni, a mathematician, and Jaros\l awa ne\'e Kowalow. After finishing The  St. Hyacinthus II Gymnasium and Lyceum in Krak\'ow he studied mathematics at the Jagiellonian University. He finished his studies obtaining the master's degree in 1939. During the occupation he took part in clandestine teaching as a member of the Underground State Examination Committee for the maturity exams. After the WWII he taught in gymnasia and lycea in Krak\'ow. In the years 1951-54 he taught mathematics at the Preparatory Study of the Jagiellonian University. From 1952 to his retirement in 1972 he taught in the A. Witkowski  V Lyceum, where he also held classes for gifted youth.

His pupils remembered him fondly (\cite{Ma71}). \textit{``Not too young, stooping, of short stature, he had all characteristics of a mathematics teacher, and  was one indeed. (...) Mathematics, the nightmare and curse of generations of pupils, was perceived by us differently. To be sure, the majority in the class was far from admiration for ``the queen of sciences" and from  [the possession of]  deeper mathematical knowledge, but  at least it was not afraid of the lessons. (...) Those who liked mathematics and had no problems with it, learned a lot at school."} 
Among Bielak's pupils were Zdzis\l aw Opial, Andrzej Lasota,  W\l odzimierz Mlak, and Antoni Leon Dawidowicz, later distinguished mathematicians and university professors. Bielak was also active in the Polish Mathematical Society. He died on February 3, 1991, in Krak\'ow.\\

\subsection{Franciszek Bierski (1912-2002)}
\begin{wrapfigure}{l}{0.45\textwidth}

\begin{center}

\includegraphics[width=0.45\textwidth]{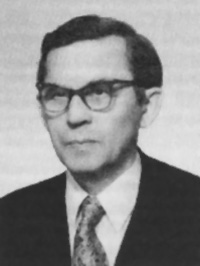}

\end{center}

\end{wrapfigure}
Born in Warszowice \'Sl\c askie, he graduated from the Jagiellonian University in 1936. He taught mathematics and physics in the gymnasium and lyceum in Piekary \'Sl\c askie. After World War II he was employed by the Academy of Mining and Metalurgy in Krak\'ow. In 1959 he got PhD in mathematics from the Jagiellonian University under the direction of Franciszek Leja. He chaired the Laboratory of Mathematical Analysis at the Academy of Mining and Metalurgy in the years 1970-83, and in 1974-83 he was the deputy director, then the director, of the Institute of Applied Mathematics. Zentralblatt f\"ur Mathematik lists 9 research publications authored or coauthored by Bierski, mostly in the field of analytic functions. He also wrote and published several academic textbooks. He died in Krak\'ow.\\

\subsection{Roman Leitner (1914-2008)} \begin{wrapfigure}{l}{0.3\textwidth} \begin{center}
\includegraphics[width=0.3\textwidth]{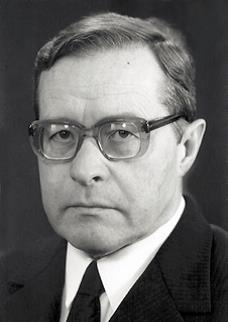}
\end{center} \end{wrapfigure} Born in 1914 in Radziech\'ow near Lw\'ow, he passed his maturity examination in Jas\l o in 1932. Then he studied
 mathematics at the Jagiellonian University, obtaining the master's degree and the high school teacher's diploma in 1937. He also studied physics (\cite{Sre07}). He was employed by the III State Gymnasium in Krak\'ow and worked voluntarily as an unpaid assistant at the Jagiellonian University, in Stanis\l aw Zaremba's chair. In September 1939 he was in Lw\'ow. 
 During the occupation he was involved in the clandestine teaching and in the years 1943-44 had to go into hiding. After the liberation of Lublin he joined the Polish Army. He became an officer of field artillery and a lecturer in the Officers' School of Shellproof Weapons (first in Che\l m, then in Modlin). As a teacher, he was released from the army in 1946 and returned to the Jagiellonian University, where he became a senior assistant. He also gave lectures in mathematics to teachers studying for professional development at the Higher Pedagogical Study in Katowice and conducted summer professional development courses for high school teachers in Szklarska Por\c eba (1949) and Ko\l obrzeg (1950).\\   
In 1949 Leitner got his PhD under the direction of Tadeusz Wa\.zewski. 
In 1951 he was called to be a deputy head of the Chair of Mathematics of newly founded Military Academy of Technology in Warsaw. The head was Witold Pogorzelski, earlier a professor of mathematics at Warsaw University of Technology, and a holder of PhD from Jagiellonian University (received in 1919).  Along with other employees, they conducted reseach in differential and integral equations and taught mathematics for applications in modern military technology. In 1954 Leitner became a \emph{docent}. In 1957 Pogorzelski returned to Warsaw University of Technology and Leitner took over the Chair. During his term (until his retirement in 1984) the Chair introduced new courses of studies, e.g.,  extramural and supplementary courses as well as programs for foreign students. Preparatory courses for applicants to technical universities  were very popular. \\

 Leitner co-organized the Television Technical University and televised preparatory courses for applicants.  Together with Wojciech \.Zakowski (later a professor in Warsaw University of Technology) he  wrote a study guide in mathematics for applicants, reprinted many times. Part of this guide became a geometry textbook for lyceum. He also wrote other lecture notes and textbooks for students, as well as educational computer programs. In 1970 the Military Academy undertook supervision of the XXIV C. K. Norwid State Lyceum (a high school) and its personnel, among them Leitner, taught advanced classes for the students. Leitner took care to organize regular instructor training for the employees of his Chair, by visiting classes and discussing performance, organizing model lectures and recitations as well as courses and seminars in methodology of teaching. He was always thoroughly prepared for his classes and his lectures were considered beautiful.  Leitner died in 2008. (\cite{Koj08a}, \cite{Koj08b}).\\

\subsection{Tadeusz Rachwa\l  \ (1914-1992)} 

\begin{wrapfigure}{l}{0.45\textwidth}

\begin{center}

\includegraphics[width=0.45\textwidth]{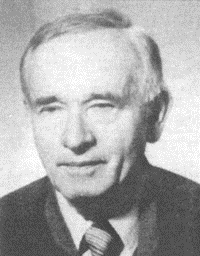}

\end{center}

\end{wrapfigure} Born on May 29, 2014 in Krak\'ow, in the family of W\l adys\l aw, a craftsman, and Maria ne\'e Kulpa. He finished A. Witkowski gymnasium in Krak\'ow in 1932 and studied at the Jagiellonian University. He credited his high school in instilling in him good studying habits in mathematics. In 1936 he obtained   the master of philoshophy degree in the field of mathematics. In October 1936 he started working in the Chair of Mathematics of the Academy of Mining in Krak\'ow, also taking up the course of studies at the Faculty of Mining. Here is how he remembered his choice of a career (\cite{Ra01}):\\

\textit{``My dreams as a graduate were partially fulfilled. I  set my sights on a different course of technical studies. I knew that I would need mathematics, so I applied myself to it with special care. But first of all I related it with my hobbies, which were painting and drawing. I wanted to take up architecture. But my fate directed me in such a way that I stayed in Krak\'ow and took up mathematics. Because later my classmates Litwiniszyn and Wojtanowicz encouraged me to study at the Academy of Mining, I decided on those studies."}\\

During the occupation Rachwa\l \  worked as a measuring technician at the Hydrological Subdivision in the region of Jas\l o, and then as an accountant in the ``Dezet" company in Krak\'ow. He also fought in the Home Army. After the war he resumed his work at the Academy of Mining, as a senior assistant. In 1950 he finished his studies at the Faculty of Mining of the Academy of Mining and Metalurgy (as the school became  known since 1949) and a year later he was nominated for the post of a deputy professor and the head of the Chair of Descriptive Geometry. In 1955, on the basis of the thesis ``A study of the order of tangency of a regular curve with a strictly tangent ball" he obtained the degree of candidate of sciences (which was at that time conferred instead of PhD). The thesis was supervised by Stanis\l aw Go\l \c ab. Rachwa\l\   got his habilitation in 1962 (at the Technical University of Krak\'ow) on the basis of the thesis ``On a certain mapping of one-sheeted hyperboloid onto a plane".  In 1971 he received the title of an extraordinary professor. He organized the Institute of Mathematics of the Academy of Mining and Metalurgy and was its first head. Until his retirement in 1984, he headed the Laboratory of Descriptive Geometry in the Institute. He initiated long-term collaboration between the Institute and the lignite mines in Turosz\'ow and Konin and the sulphur mine in Grzyb\'ow. He published about 30 research works in differential and descriptive geometry and in applications of mathematics to mining, as well as 8 textbooks and sets of lecture notes in descriptive geometry (the 2-volume text \cite{Ra73} had many editions). He supervised 7 doctorates and 4 habilitations. He died in Krak\'ow on April 18, 1992.\\

\textbf{Acknowledgments:} \\
The drawing of Antoni Bielak  made by K. Ma\-\l a\-chow\-ski was reproduced from \cite{Ma71}. The recordings of Fr. Bernard Turowicz were made available by Dr. Zofia Pawlikowska-Bro\.zek. The photos come from the archives of Z. Pawlikowska – Brożek and S. Domoradzki. They were obtained from private individuals or UJ Archives.  
The first author was partially supported by the Centre for Innovation and Transfer of Natural Sciences and Engineering Knowledge (University of Rzesz\'ow). 

\section{APPENDIX: Master's theses in mathematics prepared in Krak\'ow between 1928 and 1936}

Translated by the second author from \cite{KrUJ2334} and \cite{KrUJ3436}. Some misprints of the Polish original titles were corrected. The spelling of proper names was preserved.\\

Up to 1930/31:\\

\begin{enumerate}

\item Stanis\l aw Turski: \textit{Determining the magnitudes of quaternionic precessions by the Newcomb constants}

\item Andrzej Turowicz: \textit{On an application of iteration to solving differential and integral equations}

\item Gerson Gottesfeld: \textit{Fundamental properties of linear congruences} 

\item Micha\l \ Seidmann: \textit{Definitions of lines of curvature}

\item Etla Horn\'owna: \textit{A short outline of tangential transformations in the plane from the viewpoint of Lie}

\item Stanis\l aw Krystyn Zaremba: \textit{Differential equations and tangential transformations in the projective plane}

\item Rev. J\'ozef St\c epie\'n: \textit{Kinematic method in the theory of surfaces}

\item Jan Skar\.ze\'nski: \textit{On Prof. {\L}ukasiewicz's theory of deduction}

\item Regina Hausmann\'owna: \textit{Principles of algebra and analysis of tensors}

\item Lija Jankielowska: \textit{On inflexibility of elliptic surfaces}

\item Izrael Brumberg: \textit{A few fundamental theorems in the theory of minimal surfaces}

\item Aleksander Or\l owski: \textit{Curves of constant width}

\item Roman Dniestrza\'nski: \textit{Vectorial method in the theory of surfaces}

\item J\'ozef Steczko: \textit{The notion and properties of parallelism in a Riemann space}

\end{enumerate}

1930/31\\

\begin{enumerate}

\item Chaim Wasserfall: \textit{ Surfaces of constant curvature}

\item Wolf Kestenblatt: \textit{From the theory of analytic continuation (Mittag-Leffler star)}

\item Anna Zofja Czarkowska: \textit{Fundamental theorems in the theory of conformal transformations of planar domains}

\item Hersz H\"andel: \textit{The fundamental theorem on geodesic curves on a surface}

\item Stanis\l aw Malecki: \textit{On birational transformations}

\end{enumerate}

1931/32\\

\begin{enumerate}

\item Adam Bielecki: \textit{On integral representation of  surfaces and curves by implicit functions}

\item W\l adys\l aw Hetper: \textit{Abel-Laplace integral equations}

\item Szymon Berg: \textit{Fundamental theorems by Brill-Noether}

\item Karol Kozie\l \ : \textit{A mathematical formulation of a problem in the theory of refraction}

\item Sr. Prezepja Wilczewska: \textit{The theory of general complex numbers as an application of the theory of Lie groups}

\item Emma Epstein\'owna: \textit{Singular points of an analytic function given by a Taylor series on the circle of convergence of this series}

\item Antoni Nykli\'nski: \textit{Fundamental properties of Weingarten surfaces}

\item Klara Goldstoff\'owna: \textit{The form of a homogeneous differential equation of Fuchs type}

\item Janina Martini: \textit{Asymptotic solutions of systems of differential equations}

\item Stefan Piotrowski: \textit{Parametric from of differential equations in partial derivatives of order one}

\item Danuta Gierulanka: \textit{On periodic ordinary integral in a real variable}

\item Florjan Szozda: \textit{Fr\'echet's natural parameters}

\item Rozalja Nord\'owna: \textit{Riemann's method in  linear partial differential equations of order two, hyperbolic}

\item Emil Reznik: \textit{Set theoretic foundations of expandability of functions in the series of Bessel functions}

\end{enumerate}

1932/33\\

\begin{enumerate}

\item Bronis\l aw Czerwi\'nski: \textit{Upper and lower Perron integrals and the question of uniqueness of solutions of a system of differential equations}

\item Juljusz Keh: \textit{Bessel differential equation and main properties of Bessel and Hankel functions}

\item Zygmunt Sejud: \textit{Frenet's formulas for an $n$-dimensional Riemann space}

\item Janina Peraus\'owna: \textit{An estimate of the domain of existence of an integral of a linear nonhomogeneous partial equation of order one}

\item Rozalja Gans\'owna: \textit{Foundations of the theory of equivalence of planar figures}

\item J\'ozef Hetper: \textit{Stokes' theorem from the topological viewpoint}

\item Rudolf Wolf: \textit{Brouwer's theorem on invariance of the number of dimensions}

\item Leopold Haller: \textit{The length of a set lying on a rectifiable continuum in relation with the counting function}

\item Stefan Sedlak: \textit{On the notion of invariant}

\item Irena Wilkoszowa: \textit{Convex functions and the functional equation $f(x+1)=xf(x)$}

\item Czes\l aw Kluczny: \textit{(Essentially) two-parameter family of solutions of a differential equation $F\big ( x,y,z,\frac{\partial z}{\partial x}
\frac{\partial z}{\partial y}\big )$}

\item W\l adys\l aw Misiaszek: \textit{The role of equations attached in the reduction of linear and homogeneous differential equations}

\item Aron Teitelbaum: \textit{Tangential transformations in relation with differential equations in partial derivatives of order one}

\item Eugenjusz Ziemba: \textit{Finite tangential transformations}

\item Wincenty \L abuz: \textit{Boundary problems of an ordinary differential equation of order two}

\item Pawe\l \ Szabatowski: \textit{Contingent and paratingent}

\item Helena Mandelbaum\'owna: \textit{On the behavior of an analytic function on the boundary of the disk of convergence}

\end{enumerate}

1934/35\\

\begin{enumerate}

\item Jadwiga R\"attig: \textit{Theory of Dirichlet series}

\item Maria Kostka: \textit{Fundamental properties of regular closed spatial curves}

\item Danuta Stach\'orska: \textit{Some sufficient conditions for integral existence of an inverse transformation to a transformation of class $C$}

\item Mieczys\l aw Warcho\l   : \textit{A catalog of principles of geometry of Riemannian spaces}

\item Sydonia Kleiner\'owna: \textit{Principles of the theory of analytic sets}

\item Helena Gelber\'owna: \textit{Transformations of so-called euclidean motions and symmetries in the plane and their properties}

\item W\l adys\l aw Skrzypek: \textit{On systems of completely integrable differential equations}

\item Maria Holcherg: \textit{A boundary problem for differential equations dependent on a parameter}

\item Stanis\l aw K\c adzielawa: \textit{Malmsten's method of seeking the integrating factor for differential equations}

\item Jan Angress: \textit{Solution of a completely integrable system of differential equations by Mayer's method}

\item Antoni Bulanda: \textit{On complete extension of of functional operators}

\item Franciszek Ryszka: \textit{Main principles of the theory of automorphic functions}

\end{enumerate}

1935/36

\begin{enumerate}

\item Franciszek Bierski: \textit{Projective geometry in two-dimensional complex space}

\item Antoni Bulanda: \textit{Maximal extensions of Hermite operators in a Hilbert space}

\item Jadwiga Dymnicka: \textit{Surfaces with constant Gaussian curvature}

\item Kazimierz Gurgul: \textit{Whether a function of a complex variable corresponding to minimal surfaces according to Weierstrass' formula must be analytic}

\item Wac\l aw Juszczyk: \textit{Transformations of euclidean motions and symmetries in the plane and their properties}

\item J\'ozef Janikowski: \textit{On transformations of a differential equations in a neighborhood of a singular point}

\item J\'ozefa Konarska: \textit{On fields of rays and differential equations of order I} 

\item Karol Ka\l u\.za: \textit{On approximation of functions of one real variable}

\item Jerzy Klimonda: \textit{On a special class of natural equations for a surface in $R_3$}

\item Maria Kostka: \textit{Fundamental properties of regular closed spatial curves}

\item J\'ozef Kwieci\'nski: \textit{Development of the theory of curves in four-dimensional euclidean space $R_4$}

\item Tadeusz Kami\'nski: \textit{Fundamental theorems in the theory of tangentail transformations}

\item Jerzy \L omnicki: \textit{Implicit functions in the domain of complex variables} 

\item Franciszek Ryszka: \textit{Main outlines of the theory of automorphic functions in one variable}

\item Jadwiga R\"attig: \textit{Theory of Dirichlet series}

\item Zofia Stock\'owna: \textit{An outline of the theory of additive set functions}

\item Herbert Welke: \textit{Extremal properties of the circle and the ball}

\item Jakub Zar\c eba: \textit{Properties of the group of projective transformations of the plane}

\item Jerzy Kuczy\'nski: \textit{Tentative proofs in Prof. Wilkosz's tribe logic taking into account cardinal set theory}

\item Jan Karafia\l : \textit{On mechanical integration of certain differential equations}

\item J\'ozef Lesikiewicz: \textit{On change of variables in certain differential equations}

\end{enumerate}


\end{document}